\documentstyle[12pt]{article}
\font\tenmsb=msbm10
\font\sevenmsb=msbm7
\font\fivemsb=msbm5
\newfam\msbfam
  \textfont\msbfam=\tenmsb
  \scriptfont\msbfam=\sevenmsb
  \scriptscriptfont\msbfam=\fivemsb
\def\Bbb{\fam\msbfam\tenmsb}

\newcommand{\sd}{\mbox{${\cal{S}}$}}
\newcommand{\cP}{\mbox{${\cal{P}}$}}

\def\C{{\Bbb C}}

\def\ci{{\cal I}}
\def\O{{\cal O}}
\def\si{\sigma }

\def\z{\zeta }

\def\la{\lambda }
\def\La{\Lambda }
\def\ga{\gamma }

\def\a{\alpha }
\def\b{\beta }

\def\om{\omega }

\def\p{\phi }
\def\r{\rho }
\def\th{\theta }
\def\Th{\Theta }
\def\bx{{\tilde x}}
\def\bsi{{\tilde \sigma}}
\def\bL{{\tilde L}}
\def\ra{\rightarrow}
\def\lra{\longrightarrow}
\def\lwp{\mbox{\Large $\wp$}}

\begin{document}
\baselineskip=16.5pt

\title{Projective structures on a Riemann surface\footnote{To appear
in International Mathematics Research Notices}}

\author{Indranil Biswas ~~ and ~~ A. K. Raina}

\date{}

\maketitle

\section{Introduction}

A {\it projective structure} (also called a {\it projective
connection}) on a Riemann surface is an equivalence class of
coverings by holomorphic coordinate charts such
that the transition functions are all M\"obius transformations.
There are several equivalent notions of a projective
structure \cite{D}, \cite{G}.

For a compact Riemann surface $X$ of any genus $g$, let $L$
denote the line bundle $K_{X\times X}\otimes {\O}_{X\times
X}(2\Delta)$ on $X\times X$, where $K_{X\times X}$ is the
canonical bundle of $X\times X$ and $\Delta$ is the diagonal
divisor. This line bundle $L$ is trivialisable over a Zariski
open neighborhood of $\Delta$ and has a {\it canonical
trivialisation} over the nonreduced divisor $2\Delta$.

Our main result [Theorem 3.2] is that the space
of projective structures on $X$ is canonically identified with
the space of all trivialisations of $L$ over $3\Delta$ which
restrict to the canonical trivialisation
of $L$ over $2\Delta$ mentioned above.

In (\cite{D}, page 31, Definition 5.6 bis) Deligne gave another
definition of a projective structure (what he calls
``forme infinit\'esimale''). We give a direct
identification of this definition with our definition of a
projective structure [Theorem 4.2].

In Section 5, which is independent of the rest of the
paper, we describe briefly the
origin of this work in the study of the
so-called ``Sugawara form" of the energy-momentum tensor in a 
conformal quantum field theory.

\section{Trivialisability of the line bundle $L$}

Let $X$ be a compact connected Riemann surface, equivalently,
a smooth connected projective curve over $\C$, of genus $g$.
We denote by $S$ the complex surface $X\times X$, by $\Delta$
the diagonal divisor of $S$, and by 
$K_S$ the canonical bundle of $S$. Thus
$K_S  =  p^*_1K_X\otimes p^*_2K_X$, where $p_i$ ($i=1,2$) 
is the projection of $X\times X$ onto the $i$-th factor
and $K_X$ is the canonical bundle of $X$.

Let $\si$ be the involution of $S$ defined by $(x,y) \longmapsto
(y,x)$, of which $\Delta$ is the fixed point set. We note that $\si$
has a canonical lift $\bsi$ to $L=K_S\otimes {\O}_S(2\Delta)$; in
other words, $\bsi$ is an
isomorphism between $L$ and $\si^*(L)$ with $\bsi\circ \bsi$ being
the identity isomorphism.

The aim of this section is to establish the 
following theorem:

\noindent {\bf Theorem 2.1.}~ {\it The line bundle
$L :=  K_S\otimes {\O}_S(2\Delta)$
on $S$ is (a), trivialisable on every infinitesimal neighborhood
$n\Delta$ of $\Delta$ in $S$ and (b), has a canonical trivialisation
on the first infinitesimal neighborhood $2\Delta$, which is the
unique trivialisation of $L$ on $2\Delta$ invariant under the
action of $\bsi$ and coinciding with the canonical trivialisation
of $L$ on $\Delta$.}
  
We denote by $J^d$ the component of the Picard group of $X$
consisting of all line bundles of degree $d$ and by ${\O}(\Th)$
the line bundle on $J^{g-1}$ corresponding to the theta divisor,
viz.  the reduced theta divisor on $J^{g-1}$ defined by the
subset $\{\xi \in J^{g-1}\vert H^0(X , \xi) \neq 0\}$ (when $g
=0$, it is the zero divisor). Let $\th$ denote the natural
section of the line bundle ${\O}(\Th)$ on $J^{g-1}$ given by the
constant function $1$; it vanishes precisely on the theta
divisor.  For $\xi\in J^{g-1}$ let $\xi^*\Th$ denote the divisor
$\{\zeta\otimes
\xi^{-1}\mid \zeta\in\Th\}\subset J^0$. We now recall \cite{NR} that
the linear equivalence class of $\xi^*\Th +
(K\otimes\xi^{-1})^*\Th$ on $J^0$ is independent of $\xi\in
J^{g-1}$ and defines {\it canonically} a line bundle on $J^0$,
which we denote by ${\O}(2\Th_0)$.

We require the following property of the line bundle $L$ in the
proof of Theorem 2.1:

\noindent {\bf Lemma 2.2.}~ {\it Let $\phi~: ~~S~\lra ~J^0$ be
the morphism defined by $(x,y)\longmapsto {\cal O}_X(x-y)$. Then
$$
L~~=~~\phi^*{\O}(2\Theta_0)\leqno{(2.3)}
$$}

\noindent {\bf Proof.} Clearly we can write
$$
L ~~ =~~ {\cal M}_{\a}\otimes {\si}^*{\cal M}_{\a} \leqno{(2.4)}
$$
where, for $\a\in J^{g-1}$, the line bundle ${\cal M}_{\a}$ on $S$
is defined as follows:
$$
{\cal M}_{\a} ~ := ~ p^*_1(K_X\otimes {\a}^{-1})
\otimes p^*_2({\a})\otimes {\O}_S(\Delta) \leqno{(2.5)}
$$
As shown in \cite{R1}, ${\cal M}_\a$ is 
isomorphic to $\p_\a ^*{\O}(\Th)$, where
$$
{\p}_{\a} ~ : ~~ S ~ \lra ~ J^{g-1} \leqno{(2.6)}
$$
is the morphism defined by $(x,y) \longmapsto \a\otimes {\O}_X(x-y)$.
Theorem 2.2 of \cite{BR} is, however, preferable, since it
gives, in this special case, a {\it natural} isomorphism between
${\p}^*_{\a}{\O}(\Th)$ and ${\cal M}_{\a}\otimes {\z}_{\a}$, where
${\z}_{\a}$ denotes the trivial line bundle on $S$
with fiber ${\Th}_{\a}$, the fiber of ${\O}(\Th)$ at the point $\a$.
(For any $\a$ outside the theta divisor, the nonzero vector
${\th}(\a) \in {\Th}_{\a}$ identifies ${\z}_{\a}$ with
the trivial line bundle.)
Using this in (2.4), we see that Lemma 2.2 follows immediately
from the definition of ${\O}(2\Th_0)$.$\hfill{\Box}$

\noindent {\it Proof of part (a) of Theorem 2.1.}~ Simply observe
that the image ${\p}(\Delta) = 0\in J^0$. In fact, since
${\O}(2\Th_0)$ has no base points, $L$ has a global section
which is nowhere zero on $\Delta$.$\hfill{\Box}$

\medskip
\noindent {\bf Corollary 2.7.}~ {\it Let $\a \in J^{g-1}\setminus\Th$.
Then the section
$$
{\om}_{\a} ~~ = ~~ {\p}^*_{\a}\th \otimes
({\p}_{\a}\circ \si)^*\th ~~ \in ~~ H^0(S,L)
$$
is 1 at any point of the diagonal $\Delta$. In particular,
this section gives a trivialisation of $L$ over some Zariski open
neighborhood of $\Delta$.
The existence of ${\om}_{\a}$ implies that
$$
\dim H^0(S,L)\, \geq \, \dim H^0(S, K_S) +1 \, = \,
\dim H^0(X,K_X)^{\otimes 2} +1 \, = \, g^2+1\leqno{(2.8)}
$$}

\noindent {\bf Proof.}~ Using the natural trivialisation of ${\O}(\Th)$ 
outside the theta divisor given by the section $\th$ and the above
identification of ${\cal M}_{\a}$ with
the pullback of ${\O}(\Th)$, we have a
trivialisation of ${\cal M}_{\a}$ over
some Zariski open neighborhood of
$\Delta$. This gives a trivialisation of ${\si}^*{\cal M}_{\a}$ over
some Zariski open neighborhood of $\Delta$. Now the equality (2.4)
completes the proof.$\hfill{\Box}$
\medskip

\noindent {\bf Notation:}~ For $n\geq 1$, we shall denote the
restriction of $L$ to the divisor $n\Delta$ (the $(n-1)$-th
order infinitesimal neighborhood of $\Delta$) by $L\mid
n\Delta$.
\medskip

\noindent {\it Proof of part (b) of Theorem 2.1.}~ Now $L\mid 
\Delta={\O}_{\Delta}$ and ${\O}_{\Delta}$ has a one-dimensional
space of sections invariant under the action induced by the 
involution $\si $ on $S$. Hence $L$ has a
canonical trivialisation on $\Delta$ defined by the section ``1".

The situation on $2\Delta$ is more complicated. We know
that $\om_{\a}$ in Corollary 2.7, which defines a trivialisation
of $L$ on $2\Delta$, is symmetric under $\si$. The claim that
$L$ has a {\it canonical}
trivialisation on $2\Delta$ will then follow from the following lemma:

\noindent {\bf Lemma 2.9.}~ {\it The restriction of $L$ to $ 2\Delta$ has
a one-dimensional space of sections invariant under the action induced by
the involution $\si : (x,y)\mapsto (y,x)$ on $S$.}

\noindent {\bf Proof.}~ Consider the exact sequence
$$
0~\lra K_{\Delta}~\lra ~L\mid 2\Delta~\lra ~{\O}_{\Delta}~\lra 0
$$
where we have made use of the canonical trivialisation of $L$ on $\Delta$.
Now note that the global sections form a short exact sequence. Observe
that the natural invariant section ``1" of ${\O}_{\Delta}$ lifts
(by averaging over $\bsi$) to an invariant section of $L\mid
2\Delta$, so that the dimension of the space of invariant
sections of the latter is at least one. On the other hand,
$\bsi$ operates on $H^0(K_{\Delta})$ as {\it -Id}. Indeed, the
tangent space at $(x,x)\in\Delta$ is $T_xX \oplus T_xX$, and it
is the direct sum of the subspace spanned by $(v_x,v_x)$ with
the subspace spanned by $(v_x,-v_x)$, where $v_x$ is a nonzero
vector in $T_xX$. The former are invariant under $\bsi$ and
belong to the tangent bundle of $\Delta$, while the latter are
anti-invariant under $\bsi$ and belong to the normal bundle of
$\Delta$. Now, since $K_{\Delta}$ is the conormal bundle to
$\Delta$, the involution $\bsi$ operates as {\it -Id} on
$H^0(K_{\Delta})$. Thus we conclude that $K_{\Delta}$ has no
nonzero section which is invariant under $\bsi$. This proves the
lemma and also completes the proof of Theorem 2.1.

What is happening is that, under the quotient map $q:S\ra
S/\si$, the line bundle $L$ descends to $\bL$ on $S/\si$, since
$\bsi$ acts trivially on the fibers of $L$ at each point of the
fixed point set $\Delta$ of $\si$. The trivialisation of $L$
over $\Delta$ induces a trivialisation of $\bL$ over
$\Delta/\si$. Since the scheme-theoretic inverse image
$q^{-1}(\Delta /\si)$ is $2\Delta$ and $q^*\bL=L$, the
trivialisation of $\bL$ over $\Delta/\si$ induces a
trivialisation of $L$ over $2\Delta$.$\hfill{\Box}$

It is useful to have an alternative view of the canonical
trivialisation of $L$ on $2\Delta$:

{\bf Proposition 2.10.}~{\it The canonical trivialisation of
$L$ on $2\Delta$ is given by the unique section of $L\mid 2\Delta$
which restricts to the canonical trivialisation on $\Delta$
and lifts to a global section of $L$.}

The proof rests on the following lemma, which shows that the
inequality (2.8) is actually an equality.

\medskip
\noindent {\bf Lemma 2.11.}~~~~~ $\dim H^0(S,L) ~ = ~ g^2+1$.
\medskip

\noindent {\bf Proof.}~ In view of (2.8), we merely have to establish
the upper bound. Indeed, tensoring the following exact
sequence of sheaves on $S$
$$
0 ~ \lra ~ {\O}_S(-\Delta) ~ \lra ~{\O}_S ~
\lra ~ {\O}_{\Delta} 
~ \lra ~ 0 \leqno{(2.12)}
$$
by $L$ and passing to cohomology, this follows from the observation
that $K_S(\Delta)$ has $g^2$ sections. To establish the latter, tensor 
$(2.12)$ by $K_S(\Delta)$
and pass to cohomology; it then suffices to show that the injection
$H^0(S,K_S)\ra H^0(S,K_S(\Delta))$ is an isomorphism.
Taking the direct image of this short exact sequence by the
projection, $p_1$, to the first factor of $S$, gives the long
exact sequence 
$$
0 ~ \lra ~K_X\otimes \C^g ~\stackrel{\iota}{\lra} ~
p_{1*}(K_S(\Delta))~\lra ~K_X~ \lra ~\cdots
$$
of which the first three terms are locally free sheaves.
The first two terms are rank $g$ vector bundles and
hence $\iota$ must be an isomorphism, which completes the proof. 
$\hfill{\Box}$
\medskip

\noindent {\it Proof of Proposition 2.10.}~ From the exact sequence
$$
0~\lra ~K_S~\lra ~L~\lra ~ L\mid 2\Delta ~\lra~0
$$
and the fact that $L$ has only $g^2 + 1$ sections, we conclude
that the space of sections of $L$ has a one-dimensional image
in $L\mid 2\Delta$.$\hfill {\Box}$

\section{Projective structures and the line bundle $L$}

We will recall the definition of
a {\it projective structure on a Riemann surface subordinate to the
complex structure}. This is defined (see \cite{G} page 167) to
be a holomorphic
coordinate covering, $\{U_i,z_i\}_{i \in I}$, of $X$
such that for any pair $i, j \in I$, the holomorphic transition
function $f_{i,j}$ (defined by $z_i = f_{i,j}(z_j)$)
is a M\"obius transformation, i.e., a function of the form
$$
z ~~ \longmapsto ~~ {{az+b}\over {cz+d}} \leqno{(3.1)}
$$
where $a,b,c,d \in \C$ with $ad-bc = 1$. The space ${\lwp}$ of all
projective structures on $X$ subordinate to the
complex structure is an affine space for the complex vector space
$H^0(X, K^2_X)$ (\cite{G} page 172).

The main result of this section is the following theorem:

\noindent {\bf Theorem 3.2.}~{\it Let ${\cal Q}$ denote the
space of all trivialisations of $L\mid 3\Delta$, which, on
restriction to $2\Delta$, give the canonical trivialisation of
$L\mid 2\Delta$. Then ${\cal Q}$ is an affine space for the
vector space $H^0(X, K_X^2)$, which is canonically isomorphic to
the affine space $\lwp$ of projective structures on $X$.}

{\bf Proof.}  The obvious exact sequence
$$
0 ~ \lra ~ K^2_X ~ \lra ~ L\mid 3\Delta ~
\lra~ L\mid 2\Delta ~ \lra ~ 0 \leqno{(3.3)}
$$
shows that ${\cal Q}$ is an affine space for the vector space
$H^0(X,K^2_X)$.

We shall now construct a map from ${\lwp}$
to ${\cal Q}$.

Let $M = CP^1\times CP^1$, and consider the {\it trivial} line bundle
$L_M  ~ := ~ K_M \otimes {\O}_M(2{\Delta}_M)$
on $M$, where ${\Delta}_M$ is the diagonal on $M$. Let
$$
s ~\in ~ H^0(M,L_M) \leqno{(3.4)}
$$
be the trivialisation of $L_M$ whose restriction to ${\Delta}_M$
coincides with the canonical trivialisation given by Theorem 
2.1(b). The group of all automorphisms of $CP^1$, namely
${\rm Aut}(CP^1)$,
acts naturally on $M$ by the diagonal action; this action lifts
to $L_M$. The section $s$ in (3.4) is evidently invariant
under the induced action of ${\rm Aut}(CP^1)$ on $H^0(M, L_M)$.
This
invariance property of the section $s$ immediately implies that
if we have a projective structure on $X$, the section $s$ induces a
trivialisation of $L$ on some analytic open neighborhood of the
diagonal $\Delta$. Now, restricting this trivialisation of $L$ to
$3\Delta$ we get an element in $\cal Q$. This gives the required
map
$$
F~: ~ \lwp ~ \longrightarrow ~ {\cal Q} \leqno{(3.5)}
$$
The above construction of the map $F$ has been motivated by
\cite{Bi}.

The proof of Theorem 3.2 is now completed by the following lemma
which describes how the map $F$ relates the affine structures on
$\lwp$ and ${\cal Q}$.

\medskip
\noindent {\bf Lemma 3.6.}~ {\it For any $\ci\in {\lwp}$ and
$\ga \in H^0(X, K^2_X)$, the following equality holds.
$$
F(\ci +{\ga}) ~~ =~~ F(\ci ) ~+ ~{\ga\over 6}
$$}
\medskip

\noindent {\bf Proof.}~  Let us first recall
how the affine $H^0(X, K^2_X)$ structure of ${\lwp}$ is defined
(\cite{G} page 170, Theorem 19). We start by recalling
the definition of the {\it Schwarzian derivative},
denoted by $\sd$, which is the differential operator:
$$
\sd (f)(z)  ~:= ~ {{2f'(z)f'''(z) - 3(f''(z))^2}
\over {2(f'(z))^2}}
$$
defined over $\C$.

Take any $\ci = 
\{U_i,z_i\}_{i\in I} \in {\lwp}$ and $\ga \in H^0(X,K^2_X)$.
On each $U_i$ there is a 
holomorphic function $h_i$ such that $\ga = h_i dz_i\otimes dz_i$.
For $i \in I$, let $w_i$ be a holomorphic function on $z_i(U_i)$
satisfying the equation
$$
h_i ~= ~ \sd (w_i)(z_i) \leqno{(3.7)}
$$
Another function $w'_i$ satisfies the equation (3.7) if and only if
$w'_i (z_i) = \r\circ w_i(z_i)$, where $\r$ is a M\"obius
transformation. The element $\ci + \ga \in {\lwp}$ is given by $\{U_i, 
w_i\circ z_i\}_{i \in I}$. (Actually we may have to shrink each
$U_i$ a bit so that $w_i\circ z_i$ is a coordinate function.)

We require an explicit description
of the section $s$ defined in (3.4) in terms of local coordinates.
Identify
$CP^1$ with $\C \cup \{\infty\}$ and let $(z_1,z_2)$ be the natural
coordinates on $M$. In these coordinates the section $s$ can
be written as:
$$
s_z ~~:= ~~{{dz_1\wedge dz_2}\over {(z_1-z_2)^2}} \leqno{(3.8)}
$$

Let $ \ci ~ := ~ \{ U_i,z_i\}_{i\in I}$ be a projective 
structure on $X$, as before. Take a coordinate chart $(U , z)$ in $\ci$.
On the open set $U\times U \subset S$ there is a
natural coordinate function $(z_1,z_2)$ obtained from $z$. Now
$s_z$ in (3.8) gives a trivialisation of $L$ over $U\times U$.

Let $(V,y)$ be another coordinate chart in $\ci$ with
$y  = (az+b)/(cz+d)$ as in (3.1). This
implies that the following identity holds:
$$
s_z ~:= ~ {{dz_1\wedge dz_2}\over {(z_1-z_2)^2}} ~~  = ~~
{{dy_1\wedge dy_2}\over {(y_1-y_2)^2}} ~ =: ~ s_y
$$
where $(y_1,y_2)$ is the coordinate function on $V\times V$.
This equality implies that the two local sections of $L$,
viz. $s_z$ and $s_y$, coincide on the intersection
$(U\cap V) \times (U\cap V) \subset S$.

Thus various local trivialisations of $L$ of the
form $s_z$ patch together to give a trivialisation of $L$ on
some analytic open neighborhood of $\Delta$. In particular,
we get a trivialisation of $L\mid n\Delta$ for any $n$. Since
the section $s_z$ takes the value $1$ on $U\times U$ and is
invariant under the involution
$\bsi$, the trivialisation of $L\mid 2\Delta$ obtained this way
is the canonical trivialisation. Thus the trivialisation
of $L\mid 3\Delta$ is actually an element of ${\cal Q}$.
Evidently, the element of $\cal Q$ obtained in this way coincides
with $F(\ci)$, where $F$ is the map defined in (3.5).

Let $(U_i,z_i)$, $i\in I$, be a coordinate chart around $x \in X$,
with $z_i(x) = 0$.
Assume that
$$w_i ~= ~ z_i + \sum_{j=2}^{\infty}a_jz_i^j \leqno{(3.9)}
$$
is a solution of (3.7) (we may assume that $w_i$ is of this form
since we may 
compose $w_i$ with any M\"obius transformation).
Then equation (3.7)
gives
$$
h(0) ~=~ 6a_3 -6a^2_2 \leqno{(3.10)}
$$

Set $y = w_i\circ z_i$, and define $s_y$ as in (3.8). For
$\bx := (x,x) \in S$, using (3.9) we find that
$$
s_y(\bx ) ~=~ s_{z_i}(\bx ) + (a_3 -a^2_2)dz_i\otimes dz_i
$$
Comparing this with (3.10) we get that $s_y(\bx ) =  s_{z_i}(\bx )
+ \ga (x)/6$. This completes the proof of the lemma and 
also of Theorem 3.2.$\hfill{\Box}$ \medskip

\noindent {\bf Remark 3.11.}~ A consequence of Theorem 3.2
is the following alternative
definition of the Schwarzian derivative.

Let $f$ be a holomorphic function around $z_0 \in \C$ such that
$f'(z_0) \neq 0$. Then the function ${\bar f} := (f,f)$ is a
biholomorphism defined on some neighborhood, $U$, of
$(z_0,z_0) \in \C\times\C$. Consider the section $s$ defined in
(3.4). The restriction of
$$
{\hat s} ~ := ~ {\bar f}^*s ~ - ~ s
$$
to the (nonreduced) divisor $3{\Delta}_U$ is actually a local
section of $K^{2}_{\C}$ around
$z_0$. From the computation in the proof of Lemma 3.6 it follows
that $\hat s$ is actually ${\sd}(f)(dz)^{\otimes 2}/6$.
\medskip

\noindent {\bf Remark 3.12.} An interesting question, arising
naturally from Theorem 3.2, is whether an element of ${\cal Q}$
comes necessarily from a global section of $L$. Thus
let $\La \subset H^0(S,L)$ denote the affine subspace consisting
of those sections of $L$ which restrict to the canonical
trivialisation on $2\Delta$. Then $\La$ is an affine space
for the subspace $H^0(S,K_S)=H^0(X,K_X)^{\otimes 2}$.
Associating to any $s \in \La$, the corresponding trivialisation
of $L$ over $3\Delta$, we get a map from
$\La$ to ${\cal Q}$. Then from Theorem 3.2 we have a
(holomorphic) map $\la$ from $\La$ to $\lwp$, the space of all
projective structures on $X$. Our question is now whether $\la$
is surjective. Let
$$
R ~: ~   H^0(X, K_X)^{\otimes 2} ~= ~
H^0(S,K_S) ~\lra ~ H^0(\Delta , K_S{\vert}_{\Delta}) ~ = ~
H^0(X, K^2_X) \leqno{(3.13)}
$$
denote the obvious restriction
map. From Lemma 3.6 it follows that for any $s \in \La$ and $\b
\in H^0(X,K_X)^{\otimes 2}$, the equality $$
\la (s + \b ) ~ = ~ \la (s) + R(6\b) ~\in ~ \lwp \leqno{(3.14)}
$$
holds. This equality implies that $\la$ is surjective if and only if
the homomorphism $R$ in (3.13) is surjective. From M. Noether's theorem
(\cite{ACGH}, page 117) we know that if $X$ is non-hyperelliptic then
$R$ is surjective. Moreover, for elements $s,t\in \La$ to have
the same image under $\la$, we must have $s-t\in H^0(S,K_S(-\Delta))$.
A similar argument shows that in the non-hyperelliptic case $\la$
remains surjective when restricted to the subspace of $\La$ consisting
of sections symmetric under the map $\bsi$ induced from $\si:S\ra S
~((x,y)\mapsto (y,x))$. Related observations have been made by
Tyurin \cite{T}.
\medskip

\noindent{\bf Remark 3.15.}~ Let $(X_T, {\Gamma}_T) \longrightarrow T$
be a family of Riemann surfaces
with theta characteristic. This means that $\Gamma_T$ is a holomorphic
line bundle on $X_T$, the total space of the family of Riemann
surfaces, and for any $t\in T$, the restriction
${\Gamma}_t$ to the Riemann surface $X_t$ satisfies the condition that
${\Gamma}^{\otimes 2}_t = K_{X_t}$. Consider the line bundle
$$
{\cal M}_T ~ := ~ p^*_1(K_{\rm rel}\otimes {\Gamma}^*_T) \otimes
p^*_2({\Gamma}_T)\otimes {\O}(\Delta_T)
$$
on the fiber product $X_T\times_TX_T$, where $p_i$ denote the projection
to the $i$-th factor, ${\Delta}_T$ is the diagonal divisor in
the fiber product, and $K_{\rm rel}$ is the relative canonical
bundle on $X_T$. Let ${\cal M}_t$ be the line bundle on
$X_t \times X_t$ obtained by setting $\a = \Gamma_t$ in the proof of Lemma 
2.2. Clearly the restriction of ${\cal M}_T$ to $X_t\times X_t$ is
${\cal M}_t$. The natural isomorphism between
${\phi}^*_{\a}{\O}(\Theta)$ and ${\cal M}_{\a}\otimes {\z}_{\a}$
mentioned in the proof of Lemma 2.2 shows that the restriction
of ${\cal M}_T$ to the diagonal $\Delta_T$ is the trivial line bundle.
We may extend this trivialisation to some analytic neighborhood
of $\Delta_T$.
Now using the equality (2.4) for the given family of Riemann surfaces
we get a holomorphic family of trivialisations of
the restriction of $L$ to some neighborhood of the diagonal.
Using Theorem 3.2 this family of trivialisations
equips the family $X_T$ with a holomorphic family of
projective structures.

Given a family of Riemann surfaces, $X_{T'} \longrightarrow T'$, consider
the corresponding family of Riemann surfaces with theta characteristic
$$
(X_T ,{\Gamma}_T) ~\longrightarrow ~T
$$
where $p : T\longrightarrow T'$ is the finite \'etale Galois cover with
the fiber of $p$ over $t \in T'$ being the set of all theta 
characteristics on the corresponding Riemann
surface $X_t$. We earlier saw that there is a holomorphic family of 
projective structures for the family $X_T \longrightarrow T$. For $x \in 
T$ let ${\lwp}_x$ denote the projective structure on the Riemann surface
over $x$. For any $t \in T'$ consider the projective structure
on $X_t$ given by the average
$$
{1\over {\# p^{-1}(t)}}\sum_{x \in p^{-1}(t)} {\lwp}_x
$$
which is defined using the affine space structure on the space of all
projective structures on $X_t$. Using this construction we conclude that
the family of Riemann surfaces, $X_{T'}$, admits a holomorphic
family of projective structures.

\section{Relation with Deligne's definition}

We shall now recall another definition of a projective structure
given in \cite{D} (Definition 5.6 bis).

The fibers of the natural
projection, $\nu$, of the second order infinitesimal neighborhood
of the diagonal $\Delta$ (in $S$) onto $\Delta$ are isomorphic to
${\rm Spec}(R)$, where $R$ is the algebra ${\C}[\epsilon
]/{\epsilon}^3$. Let $P$ denote the
principal ${\rm Aut}({\rm Spec}(R))$ bundle on $X$ whose fiber over
$x \in X$ is the space of all isomorphisms between ${\rm Spec}(R)$
and the fiber of $\nu$ over $x$. On the other hand, ${\rm Aut}({\rm 
Spec}(R))$ is same as the group of all
automorphisms of $CP^1$ that fix the point $0 \in \C \cup
\{\infty\} = CP^1$. Let $P_{tg}$ denote the
projective bundle on $X$ associated to $P$. Since ${\rm
Aut}({\rm Spec}(R))$
fixes a point in $CP^1$, the bundle $P_{tg}$ has a natural
section which we shall denote by $\tau$. There is a natural
isomorphism between the second order infinitesimal neighborhood of 
$\Delta$ and
the second order neighborhood of the image of $\tau$ (in $P_{tg}$).

\medskip
\noindent {\bf Definition 4.1} ``Definition 5.6 bis of
\cite{D}''. A projective structure on $X$ is
an isomorphism between the third order infinitesimal neighborhood
of the diagonal $\Delta$ (in $S$) with the third order infinitesimal
neighborhood of $\tau$
(in $P_{tg}$) such that
the restriction of this isomorphism to the second order 
infinitesimal neighborhood
of $\Delta$ is the canonical isomorphism with the second order
infinitesimal neighborhood of the image of $\tau$ mentioned above.
\medskip

If $X = CP^1$, the projective line, then $P_{tg} = CP^1\times CP^1$.
Thus there is a canonical projective structure on $CP^1$
in the sense of (\cite{D}, Definition 5.6 bis) given by the
identity map of of the third order neighborhood of the diagonal
in $CP^1 \times CP^1$.

Let $\cal H$ denote the sheaf on $X$ which to any
open set, $U \subset X$, associates the space of all embeddings
of the third order infinitesimal neighborhood of the diagonal
of $U\times U$
into the restriction of $P_{tg}$ to $U$ which lift the canonical
embedding of the second order neighborhood of the diagonal of
$U\times U$.

Let us recall  \cite{D} that a
$K^2_X$-{\it torsor} is a holomorphic fiber bundle over $X$ such that
its fiber over any $x\in X$ is equipped with a free, transitive
holomorphic action of the fiber $K^2_x$. In other words, the result of
the action of a local holomorphic section of $K^2_X$ on a local 
holomorphic section of the torsor is again a local holomorphic section.

Proposition 5.8 (page 32) of \cite{D} says that the
sheaf $\cal H$ defined above is a $K^2_X$-torsor.

We shall denote the restriction of $L$ to $n\Delta$ by $L(n)$.
For an analytic open set $U$ of $X$ let ${\Delta}_U$ be the
diagonal divisor on $U\times U$. Let $\cal G$ denote the sheaf
on $X$ which to any open set
$U \subset X$ associates the space of all
trivialisations of the restriction of $L(3)$ to $3\Delta_U$
giving the canonical trivialisation on $2\Delta_U$. From (3.3)
it follows that the restriction ${\cal G} (U)$ is an affine
space for $H^0(U,K^2_U)$, where $K_U$ is the canonical bundle of
$U$. In other words, $\cal G$ is a torsor for the sheaf $K^2_X$.

Our aim in this section is to prove the following theorem:

\medskip
\noindent {\bf Theorem 4.2.}~ {\it The two
$K^2_X$-torsors on $X$, namely $\cal G$ and $\cal H$,
are canonically isomorphic.}
\medskip

Theorem 4.2 gives a natural identification of the
space ${\cal Q}$, the space of global
sections of $\cal G$, with the space of global sections of $\cal H$,
which is the space of all projective structures on $X$ in
the sense of (\cite{D}, Definition 5.6 bis).
\medskip

\noindent {\it Proof of Theorem 4.2.}~ We shall prove the
theorem by constructing a third $K^2_X$-torsor, $\cal T$, on $X$
and identifying both $\cal G$ and $\cal H$ with $\cal T$.
A reason for introducing $\cal T$ as the intermediate step is that
its construction might be of some interest.

For $n\geq 0$, let $J^n(X)$ denote the sheaf of jets of
order $n$ on $X$, which is a vector bundle on $X$ of rank $n+1$.
Define ${J}^n_0(X)$ to be the kernel of the obvious projection
of $J^n(X)$ onto $J^0(X)$.
Note that there is a canonical splitting of
the inclusion of ${J}^n_0(X)$ into $J^n(X)$ given by the constant
functions. Let $\cP (X)$ denote the subset of the total
space of ${J}^3_0(X)$ given by the
inverse image of $\{K_X -0\}$ (the
set of all nonzero vectors in the total space of $K_X$) under the
projection of ${J}^3_0(X)$ onto ${J}^1_0(X) = K_X$. The
space $\cP (X)$ admits a natural action (by composition of
functions) of the group $M(0)$, the
isotropy group of $0 \in \C$ for the M\"obius group action on
$CP^1$. The action
of $M(0)$ on $\cP(X)$ is free, since the only M\"obius
transformation of $CP^1$, which acts as the
identity map on the second order neighborhood of a point, is actually
the identity transformation (\cite{D}, page 29).

Let $\cal T$
denote the quotient of $\cP (X)$ by $M(0)$. A projective 
structure on $X$ gives maps from neighborhoods of points of $X$
into $CP^1$ which differ only by a M\"obius transformation, and
hence gives a section of the obvious projection 
of $\cal T$ onto $X$. An identification between the 
space of all sections of $\cal T$ and $\lwp$, the space of projective
structures on $X$, is obtained in this way.

We shall now give a $K^2_X$-torsor structure on $\cal T$.
Let $f \in {\cP}(X)$ be an element over $x \in X$, and let $v \in 
(K^2_X)_x$ be an element of the fiber of $K^2_X$ over $x$. Let
${\bar f}$ be a function defined around $x$ which represents $f$.
Since $d{\bar f}(x) \neq 0$ (by the definition of ${\cP}(X)$),
there is a number $\la \in \C$ such that $v = \la .d{\bar f}(x)\otimes
d{\bar f}(x)$. Consider the function
$$
{\bar f}_{\la}~ := ~  {\bar f} ~+ ~ \la . {\bar f}^3 \leqno{(4.3)}
$$
defined around $x$. The element in ${\cP}(X)$ over $x$
represented by the function ${\bar f}_{\la}$ clearly does not
depend upon the choice of the representative $\bar f$ of $f$.
An action of $K^2_X$ on ${\cP}(X)$ is obtained by
mapping the pair $(v ,f)$ to the element of ${\cP}(X)_x$
represented by ${\bar f}_{\la}$. This is a free (but not transitive)
action of the
abelian group scheme $K^2_X$ over $X$. This action of $K^2_X$ on
${\cP} (X)$ induces a $K^2_X$-torsor structure on the
quotient space $\cal T$ of ${\cP}(X)$. Indeed, this is a consequence
of the following fact: let $J^n_0(0)$ be the jets of order $n$ of 
functions vanishing at $0\in \C$; in this notation, the group 
$M(0)$ acts freely and transitively on
the subset of $J^2_0(0)$ consisting of all elements 
whose image in $J^1_0(0)$ is nonzero. (If $i$ denotes the isomorphism
from the space of all sections of $\cal T$ to $\lwp$, then
$i(A +\ga) = i(A)+ 6\ga$ for any $\ga \in H^0(X,K^2_X)$.)

Theorem 4.2 is a consequence of the assertion that both $\cal G$
and $\cal H$ coincide with this $K^2_X$-torsor $\cal T$. We shall first
show that $\cal G$ coincides with $\cal T$.

Take any $x \in X$, and let $f \in {\cal T}_x$ be an element of the
fiber over $x$. Let $z$ be a function defined in a
neighborhood, $U$, of $x$, that represents $f$.
Since $d{z}(x) \neq 0$, we may
assume that $z$ is a biholomorphism onto its image. Let ${\bar
z} = (z,z)$ be the biholomorphism defined on $U\times U$.
Pull back the section $s$ (defined in (3.4)) to $U\times U$ using
this map $\bar z$. Let $\hat f$ denote the local
section of $\cal G$ obtained by restricting this section to the second
order infinitesimal neighborhood of the diagonal. The 
evaluation at $x$, namely ${\hat f}(x)$, depends
only on $f$ and not on the representing function $z$. Thus we
have a map from $\cal T$ to $\cal G$ which is evidently an
isomorphism. We want to check that this isomorphism preserves
the $K^2_X$-torsor structures of $\cal T$ and $\cal G$.

Take an element $v = \la (dz)^{\otimes 2} \in K^2_x$, where $\la
\in \C$. From the
definition of the $K^2_X$-torsor structure on $\cal T$ in (4.3)
it follows that the local function $z +\la z^3$ represents the
result of the action of $v$ on $f$. Remark 3.11 says that the two
sections of $\cal G$,
represented by $z$ and $z+\la z^3$ respectively, differ by $\sd
(z + \la z^3)(0)/6$. Since
$$
\sd (z+\la z^3)(0) ~ = ~ 6\la
$$
the preservation of the $K^2_X$-torsor structures of $\cal T$
and $\cal G$ is established. 

Next we want to show that $\cal H$ coincides with $\cal T$.
Take $x$, $f$ and $z$ as above. We noted earlier that $CP^1$ has
a canonical projective structure (in the sense of \cite{D},
Definition 5.6 bis) given by the identity map of the third order
neighborhood of the diagonal. This projective structure induces
a projective structure on $U$ by the biholomorphism $z$. The
evaluation of the section (over $U$) of $\cal H$, thus obtained,
at the point $x$, does not depend upon the choice of the
representative $z$ of $f$. This gives the required
$K^2_X$-torsor structure preserving isomorphism
between $\cal T$ and $\cal H$.

As an alternative proof of Theorem 4.2 we shall give a
direct identification between $\cal G$ and $\cal H$ using
coordinate charts.

Let $(U,z)$ be a coordinate chart around $x \in X$ with
$z(x) = 0$. Using (3.8) we get a section of $\cal G$ over $U$.
We shall denote this section as $f_z$.

Since the only M\"obius transformation of $CP^1$, which acts as the
identity map on the second order neighborhood of a point, is actually
the identity map, and the group of M\"obius transformations
acts transitively on $CP^1$, there is a natural projective
structure on $CP^1$ in the sense of (\cite{D}, Definition
5.6 bis).

Since the function $z$ identifies $U$ with an open set in
$CP^1$, we get a local section of $\cal H$, which we shall
denote by $g_z$.

By mapping the section $f_z$ to $g_z$ we get an identification of
the restriction ${\cal G} (U)$ with ${\cal H} (U)$ which preserves
the torsor structures. We shall show that this identification does
not depend upon the choice of the coordinate function $z$.

Let $(V,w)$ be another coordinate chart around $x$. Thus
$$
w ~ = ~ \sum_{i=0}^{\infty}a_i z^i \leqno{(4.4)}
$$
with $a_1 \neq 0$. Let $f_w$ (resp. $g_w$) denote the local
section of $\cal G$ (resp. $\cal H$) for $(V,w)$.
It is a simple calculation using (4.4) to check that
$$
f_w(x)\, - \, f_z(x)~= ~ {{a_1a_3 - a^2_2} \over
{a^2_1}} dz\otimes dz
~ = ~ g_w(x) \, - \, g_z(x)
$$
This completes the proof of the theorem.$\hfill{\Box}$

\section{Genesis in conformal field theory}

In this section we explain how the above definition of projective 
connection in terms of trivialisations of $L$ on $3\Delta$ came out
of some investigations on a model quantum field theory on a curve
(see \cite{R1}-\cite{R3}), which give it the intuitive picture of a
{\it generalized cross ratio} on a compact Riemann surface, in the
limit when all of its arguments are made to coalesce.

The application of algebraic geometry to quantum field theory in
\cite{R1}-\cite{R3} rests on replacing the study of ``quantum
fields", which are not geometric objects, by their so-called
``$n$-point functions" which are hypothesised to be so. Thus in
\cite{R1} and
\cite{R2} we identified the ``$n$-point functions" of the defining
``quantum fields" of the model with meromorphic sections of
certain line bundles on the $n$-fold Cartesian product of the
curve $X$.  In \cite{R3} we showed how the $n$-point functions
of the ``current" $j$, which is a ``regularised product" of the
defining fields, could be computed by the use of {\it schemes}
having {\it nilpotent elements} to give a precise meaning to the
coalescing of arguments involved in the definition of the
current.

A similar regularised product of currents gives the ``energy-momentum
tensor" $T$ of the system, a fact usually expressed by
saying that $T$ is in ``Sugawara form". This is a feature 
of many conformal quantum field theories and plays an important role
in the theory of the Virasoro algebra \cite{KR}. The heuristic
expectation in (conformal) quantum field theory \cite{BPZ} is
that its ``one point function" $<T(z)>$  is a {\it projective 
connection}.

Our study of $<T(z)>$ proceeds from the
calculation of the two point function of currents $<j(z)j(w)>$ in
\cite{R3}. The salient point is the introduction of the
remarkable line bundle ${\cal A}:={\cal O}(D_{12}+D_{34}-D_{14}-D_{23})$
on $X^4:=X_1\times X_2\times X_3\times X_4$, the product of 4 copies of X,
where $D_{ij}$ denotes the divisor of $X^4$ defined by the diagonal
of $X_i$ and $X_j$. It was pointed out in \cite{R3} that the
canonical meromorphic section $1_{\cal A}$, associated with the divisor
defining ${\cal A}$, is a natural generalisation
to an arbitrary compact, connected Riemann surface of the
{\it cross ratio} of 4 points in the complex plane. It was shown in
\cite{R3} that the calculation of
$<j(z)j(w)>$ requires the trivialisability of ${\cal A}$ on 
the product scheme $Z:=2\Delta_{13}\times 2\Delta_{24}$, where 
$\Delta_{ij}$ is the diagonal of $X_i\times X_j$ and $2\Delta_{ij}$
denotes its first infinitesimal neighborhood.

\noindent {\bf Proposition 5.1}.~ {\it The line bundle
${\cal A}:={\cal O}(D_{12}+D_{34}-D_{14}-D_{23}) $
is trivialisable on $Z~ := 2\Delta_{13}\times 2\Delta_{24}$
and, moreover, if $\rho\in H^0(Z,{\cal A}\mid Z)$ denotes such
a trivialisation, then
$$
1_{\cal A}\mid Z~-~ \rho = \omega_B \leqno{(5.2)}
$$
where $\omega_B$ denotes a symmetric meromorphic section of 
$K_{X\times X}$
with double pole on the diagonal, defined by a holomorphic section of 
$K_{X\times X}(2\Delta)$ which restricts to $1$ on the diagonal.}

As pointed out in \cite{R3}, equation (5.2) can be regarded as the
precise algebro-geometric formulation of the following
well known formula
expressing the meromorphic bidifferential $\omega_B$ in terms
of the ``prime form" $E(x,y)$ (see Fay\cite{F}, eqn.(28) p.20):
$$
\omega_B(x,y)~=~\frac{\partial^2\ln E(x,y)}{\partial x\partial 
y}\leqno{(5.3)}
$$
In this way it was shown in \cite{R3} that the two point
function of currents $<j(z)j(w)>$ is a {\it symmetric
meromorphic bidifferential with a double pole on the diagonal}.
The computation of the ``one point function" $<T(z)>$ from
$<j(z)j(w)>$ now requires that the line bundle $K_{X\times
X}(2\Delta)$ should be trivialisable on the {\it second}
infinitesimal neighborhood $3\Delta$ of $\Delta$ in $X\times X$
(see \cite{R4} for further details), the validity of which
follows from results in \cite{R1}. In this way we arrive at our
proposed definition of projective connection, having started
with the generalised cross ratio and ended with the coalescing
of all of its arguments.

The fact
that a symmetric meromorphic bidifferential gives rise to a projective
connection appears to have been first observed in \cite{HS} (see also
\cite{F} p.20, following eqn.(28) cited above). The techniques
used in these references, however, do not give a {\it characterisation}
of a projective structure, as is provided by Theorem 3.2, nor an
understanding of when all projective structures arise in this way,
as is provided by Remark 3.12. Moreover, their approach cannot be 
adapted to the study of $<T(z)>$, for which we require an
algebro-geometric approach, which will make possible the study of the
{\it higher} point functions as well as other related problems. It also 
appears to be a fact that \cite{HS} and \cite{F} are
inaccessible to most geometers and so we hope that the present treatment
clarifies some of these results.

The present formulation of the concept of projective
connection was announced in several conferences and also in \cite{R4},
where the interested reader will find, in addition, a survey for
mathematicians of the papers \cite{R1}, \cite{R2} and \cite{R3}.

\medskip
\noindent {\bf Acknowledgments:} The authors are very grateful 
to Prof. M. S. Narasimhan for his useful comments. The
first named author is thankful to the Institut Fourier and the
Acad\'emie des Sciences, Paris, for their hospitality and support.
The second named author thanks the International
Centre for Theoretical Physics, Trieste, for its hospitality.
%%%%%%%%%%%%%%%%%%%%%%%%%%%%%%%%%%%%%%%%%%%%%%%%%%%%%%%%%

\noindent School of Mathematics, Tata Institute of
Fundamental Research, Homi Bhabha Road, Bombay 400005, INDIA\\
E-mail : indranil@math.tifr.res.in

\noindent Theoretical Physics Group, Tata Institute of
Fundamental Research, Homi Bhabha Road, Bombay 400005, INDIA\\
E-mail : raina@theory.tifr.res.in

\end{document}